\input amstex
\documentstyle{amsppt}
\TagsOnRight \NoRunningHeads
\magnification=\magstep1 \vsize 22 true cm \hsize 16 true cm

\topmatter
\title Conjugacy separability of some one-relator groups
\endtitle
\author D.~Tieudjo and D.~I.~Moldavanskii \endauthor


\abstract Conjugacy separability of any group of the class of
one-relator groups given by the presentation $\langle a, b;\
[a^m,b^n]=1\rangle$ ($m,n>1$) is proven.
\endabstract

\endtopmatter

\document

\centerline{ \bf Introduction}

\medskip

A group $G$ is conjugacy separable if for any two non-conjugate
elements $f$ and $g$ of $G$ there exists a homomorphism $\varphi$
of group $G$ onto some finite group $X$ such that the images
$f\varphi$ and $g\varphi$ of elements $f$ and $g$ are not
conjugate in group $X$. It is clear that any conjugacy separable
group $G$ is residually finite (i.~e., recall, for any
non-identity element $g\in G$ there exists a homomorphism
$\varphi$ of group $G$ onto some finite group $X$ such that
$g\varphi\neq 1$).

Since 1962 when G.~Baumslag and D.~Solitar [3] discovered the
first examples of non-residually finite one-relator groups a lot
of results establishing the residual finiteness of various
one-relator groups have appeared. It was also shown that some of
such groups are conjugacy separable. So, J.~Dyer [4] have proved
the conjugacy separability of any group being free product of free
groups amalgamating cycle. In the same paper she proved an
unpublished result of M.~Amstrong on conjugacy separability of any
one-relator group with nontrivial center. The conjugacy
separability of groups defined by relator of form $(a^mb^n)^t$ and
$(a^{-1}b^lab^m)^s$, where $s>1$, was proved in [6] and [1]
respectively. It should mention that up to now it is not known
whether there exists one-relator group that is residually finite
but not conjugacy separable.

In this paper we enlarge the class of conjugacy separable
one-relator groups. Namely, we prove here the

\proclaim{\indent Theorem} Any group $G_{mn}=\langle a,b;\ [a^m,
b^n]=1 \rangle$ where integers $m$ and $n$ are greater than 1, is
conjugacy separable.
\endproclaim

We note that the assertion of Theorem is also valid when $m=1$ or
$n=1$. Indeed, in this case group $G_{mn}$ is the generalized free
product of two finitely generated abelian groups with cyclic
amalgamation and its conjugacy separability follows from the
result of [4].

The residual finiteness of groups $G_{mn}$ is well known; it
follows, for example, from the result of paper [2]. Some
properties of these groups were considered in [7] where, in
particular, the description of their endomorphisms was given. The
proofs in this paper made use of presentation of group $G_{mn}$ as
amalgamated free product. The same presentation will be the
crucial tool in the proof of Theorem above and because of this, in
section 1, we recall the Solitar theorem on the conjugacy of
elements of amalgamated free products and derive, for our special
case, the criterion which is somewhat simpler. In section 2 some
properties of separability of groups $G_{mn}$ are established and
in section 3 the proof of Theorem will be completed.
\bigskip

\centerline{ \bf 1. Preliminary remarks on the conjugacy}
\centerline{ \bf in amalgamated free products}
\medskip

Let us recall some notions and properties concerned with the
construction of free product $G = (A*B;\ H)$ of groups $A$ and $B$
with amalgamated subgroup $H$.

Every element $g\in G$ can be written in a form
$$
g=x_1x_2\cdots x_r \ (r\geqslant 1)
$$
where for any $i=1, 2, \dots , r$ element $x_i$ belongs to one of
the free factors $A$ or $B$ and if $r>1$ any successive $x_i$ and
$x_{i+1}$ do not belong to the same factor $A$ or $B$ (and
therefore for any $i=1, 2, \dots , r$ element $x_i$ does not
belong to amalgamated subgroup $H$). Such form of element $g$ is
called reduced.

In general, element $g\in G$ can have different reduced forms, but
if $g=y_1y_2\cdots y_s$ is one more reduced form of $g$, then
$r=s$ and, for any $i=1, 2, \dots , r$, $x_i$ and $y_i$ belong to
the same factor $A$ or $B$. The length $l(g)$ of element $g$  is
then defined as the number $r$ of the components in a reduced form
of $g$. Element $g$ is called cyclically reduced if either
$l(g)=1$ or $l(g)>1$ and the first and the last components of its
reduced form do not belong to the same factor $A$ or $B$. If
$g=x_1x_2\cdots x_r$ is a reduced form of cyclically reduced
element $g$ then cyclic permutations of element $g$ are elements
of form $x_ix_{i+1}\cdots x_rx_1 \cdots x_{i-1}$, where $i=1,2,
\dots, r$.

If $X$ is subgroup of a group $Y$ we shall say that elements $a$
and $b$ of $Y$ are $X$-conjugate if $a=x^{-1}bx$, for some $x\in
X$.

The Solitar criterion of conjugacy of elements of amalgamated free
product of two groups (Theorem 4.6 in [5]) can be formulated as
follows.

\proclaim {\indent Proposition 1.1} Let $G = (A*B;\ H)$ be the
free product of groups $A$ and $B$ with amalgamated subgroup $H$.
Every element of $G$ is conjugate to a cyclically reduced element.
If lengths of two cyclically reduced elements are nonequal then
these elements are not conjugate in $G$. Let $f$ and $g$ be
cyclically reduced elements of $G$ such that $l(f)=l(g)$. Then
\roster \item If $f\in A$ and $f$ is not conjugate in $A$ to any
element of subgroup $H$ then $f$ and $g$ are conjugate in group
$G$ if and only if $g\in A$ and $f$ and $g$ are conjugate in $A$.
Similarly, if $f\in B$ and $f$ is not conjugate in $B$ to any
element of subgroup $H$ then $f$ and $g$ are conjugate in group
$G$ if and only if $g\in B$ and $f$ and $g$ are conjugate in $B$.
\item If $f\in H$ then $f$ and $g$ are conjugate in group $G$ if
and only if there exists a sequence of elements
$$
f=h_0, h_1, \dots h_n, h_{n+1}=g
$$
such that for any $i=0, 1, \dots , n$ $h_i\in H$ and elements
$h_i$ and $h_{i+1}$ are $A$-conjugate or $B$-conjugate. \item If
$l(f)=l(g)>1$ then $f$ and $g$ are conjugate in group $G$ if and
only if element $f$ is $H$-conjugate to some cyclic permutation of
$g$.
\endroster
\endproclaim

In the following special case assertion (2) of Proposition 1.1
became unnecessary.

\proclaim {\indent Proposition 1.2} Let $G = (A*B;\ H)$ be the
free product of groups $A$ and $B$ with amalgamated subgroup $H$.
Suppose that for any element $f$ of $A$ or $B$ and for any element
$h$ of subgroup $H$, the inclusion $f^{-1}hf\in H$ is valid if and
only if elements $f$ and $h$ commute. If $f$ and $g$ are
cyclically reduced elements of $G$ such that $l(f)=l(g)=1$ then
$f$ and $g$ are conjugate in group $G$ if and only if $f$ and $g$
belong to the same subgroup $A$ or $B$ and are conjugate in this
subgroup.
\endproclaim

In fact, let elements $f$ and $g$ be conjugate in group $G$ and
let $f$ belong to subgroup $A$, say. If $f$ is not conjugate in
$A$ to any element of subgroup $H$ then desired assertion is
implied by assertion (1) of Proposition 1.1. Otherwise, we can
assume that $f\in H$ and by assertion (2) of Proposition 1.1 there
is a sequence of elements
$$
f=h_0, h_1, \dots, h_n, h_{n+1}=g,
$$
such that for any $i=0,1, \dots, n$ $h_i\in H$ and for suitable
element $x_i$ belonging to one of the subgroups $A$ or $B$, the
equality $x_i^{-1}h_ix_i=h_{i+1}$ holds. Since for $i=0,1, \dots,
n-1$ the inclusion $f^{-1}h_{i+1}f\in H$ is valid, the hypothesis
gives $h_0=h_n$, i.~e. $x_n^{-1}fx_n=g$. This means that elements
$f$ and $g$ belong to that subgroup $A$ or $B$ which contains
element $x_n$ and also are conjugated in this subgroup.
\smallskip

Next, we describe more explicit the situation arising in assertion
(3) of Proposition 1.1.

\proclaim {\indent Proposition 1.3} Let $f=x_1x_2 \cdots x_r$ and
$g=y_1y_2 \cdots y_r$ be the reduced forms of elements $f$ and $g$
of group $G=(A*B; \ H)$ where $r>1$. Then $f$ and $g$ are
$H$-conjugate if and only if there exist elements $h_0$, $h_1$,
$h_2$, \dots, $h_r=h_0$ in subgroup $H$ such that for any $i=1,2,
\dots , r$, we have
$$
x_i=h_{i-1}^{-1} y_ih_i. \tag1
$$
\endproclaim

\demo{\indent Proof} If for some elements $h_0$, $h_1$, $h_2$,
\dots, $h_r=h_0$ of subgroup $H$ the equalities (1) hold then
$$
f=x_1x_2 \dots x_r=h_{0}^{-1}y_1h_1 \,h_{1}^{-1}y_2h_2 \, \cdots
\, h_{r-1}^{-1}y_2h_r = h_0^{-1}(y_1y_2 \cdots y_r)
h_r=h_0^{-1}gh_0.
$$

Conversely, by induction on $r$ we prove that if for some $h\in H$
the equality $f=h^{-1}gh$ holds then there exists a sequence of
elements $h_0$, $h_1$, $h_2$, \dots, $h_r=h_0$ of subgroup $H$
satisfying the equalities (1) and such that $h_0=h$.

Rewriting the equality $f=h^{-1}gh$ in the form
$$
x_r^{-1} \cdots x_2^{-1}x_1^{-1}h^{-1}y_1y_2 \cdots y_rh=1,
$$
we see that since the expression in the left part of it cannot be
reduced, elements $x_1$ and $y_1$ must be contained in the same
subgroup $A$ or $B$ and element $x_1^{-1}h^{-1}y_1$ must belong to
subgroup $H$. Denoting this element by $h_1^{-1}$, we have
$x_1=h^{-1}y_1h_1$. If $r=2$ then the equality above takes the
form $x_2^{-1}h_1^{-1}y_2h=1$, whence $x_2=h_1^{-1}y_2h$.
Therefore, setting $h_0=h_2=h$, we obtain in that case the desired
sequence.

If $r>2$ let's rewrite the equality $f=h^{-1}gh$ in the form
$$
h^{-1}y_1h_1x_2 \cdots x_r=h^{-1}y_1y_2 \cdots y_rh.
$$
This implies that if we set $f'=h^{-1}h_1x_2 \cdots x_r$ and
$g'=y_2 \cdots y_r$, then $f'=h^{-1}g'h$. Since the length of
elements $f'$ and $g'$ is equal to $r-1$, then by induction, there
exists a sequence  $h'_1=h$, $h'_2$,  \dots , $h'_r=h'_1$ of
elements of subgroup $H$, such that
$h^{-1}h_1x_2=(h'_1)^{-1}y_2h'_2$ and for any $i=2, 3, \dots , r$
$$
x_i=(h'_{i-1})^{-1}y_ih'_i.
$$
Since $x_2=h_1^{-1}y_2h'_2$, then setting $h_i=h'_i$ for $i=2, 3,
\dots , r$, we obtain the desired sequence, and induction is
complete. Proposition 1.3 is proven.
\enddemo

We conclude this section with one more property of amalgamated
free product.

\proclaim {\indent Proposition 1.4} Let $G = (A*B;\ H)$ be the
free product of groups $A$ and $B$ amalgamating subgroup $H$ where
$H$ lies in the center of each of groups $A$ and $B$. Then for any
element $f\in G$ not belonging to subgroup $A$ we have
$f^{-1}Af\cap A=H$.
\endproclaim

In fact, since subgroup $H$ coincides with the center of group $G$
(see [5], page 211), the inclusion $H\subseteq f^{-1} Af\cap A$ is
obvious. Conversely, let element $f^{-1}af$ belong to subgroup $A$
where $a\in A$. Since $f\notin A$ then either $f\in B\setminus H$
or $l(f)>1$ and (without loss of generality) the first component
of reduced form of $f$ belongs to $B$. In any case the assumption
$a\notin H$ would imply that $l(f^{-1}af)>1$. Thus, $a\in H$ and
therefore $f^{-1}af=a\in H$.
\bigskip

\centerline{ \bf 2. Some properties of groups $G_{mn}$}
\centerline{ \bf and of their certain quotients }
\medskip

In what follows, our discussion will make use of the presentation
of group $G_{mn}$ as an amalgamated free product of two groups. To
describe such presentation, let $H$ be the subgroup of group
$G_{mn}$, generated by elements $c=a^m$ and $d=b^n$. Also, let $A$
denote the subgroup of group $G_{mn}$ generated by element $a$ and
subgroup $H$, and let $B$ denote the subgroup of group $G_{mn}$,
generated by element $b$ and subgroup $H$. Then it can be
immediately verified that $H$ is the free abelian group with base
$c$, $d$, group $A$ is the free product $(\langle a \rangle * H;\
\langle a^m = c \rangle)$ of infinite cycle $\langle a \rangle$
and group $H$ with amalgamation $\langle a^m \rangle$, group $B$
is the free product $(\langle b \rangle * H;\ \langle b^n = d
\rangle)$ of infinite cycle $\langle b \rangle$ and group $H$ with
amalgamation $\langle b^n \rangle$, and group $G_{mn}$ is the free
product $(A*B;\ H)$ of groups $A$ and $B$ with amalgamation $H$.

The same decomposition is satisfiable for certain quotients of
groups $G_{mn}$. Namely, for any integer $t>1$ let $G_{mn}(t)$ be
the group with presentation
$$
\langle a, b; \ [a^m,b^n]=1, a^{mt}=b^{nt}=1\rangle
$$
and $\rho_t$ be the natural homomorphism of group $G_{mn}$ onto
$G_{mn}(t)$. Then it is easy to verify that subgroup
$H(t)=H\rho_t$ of group $G_{mn}(t)$ is isomorphic to the quotient
$H/H^t$ (where, as usually, $H^t$ consists of all elements $h^t$,
$h\in H$), subgroup $A(t)=A\rho_t$ is the amalgamated free product
$(\langle a; \ a^{mt} =1\rangle * H(t);\ a^m=cH^t)$ of cycle
$\langle a;\ a^{mt} =1\rangle$ of order $mt$ and group $H(t)$,
subgroup $B(t)=B\rho_t$ is the amalgamated free product \linebreak
$(\langle b; \ b^{nt} =1\rangle * H(t);\ b^n=dH^t)$ of cycle
$\langle b;\ b^{nt} =1\rangle$ of order $nt$ and group $H(t)$ and
group $G_{mn}(t)$ is the amalgamated free product $(A (t)*B (t);\
H(t))$.

These decompositions of groups $G_{mn}$ and $G_{mn}(t)$ are
assumed everywhere below, and such notions as free factor, reduced
form, length of element and so on will refer to them.

Let us remark, at once, that since each of groups $A(t)$ and
$B(t)$ is the amalgamated free product of two finite groups and
group $G_{mn}(t)$ is the free product of groups $A(t)$ and $B(t)$
with finite amalgamation, it follows from results of paper [4]
that for every $t$, group $G_{mn}(t)$ is conjugacy separable. So,
to prove the conjugacy separability of group $G_{mn}$ it is
enough, for any non-conjugate elements $f$ and $g$ of $G_{mn}$, to
find an integer $t$ such that elements $f\rho_t$ and $g\rho_t$ are
non-conjugate in group $G_{mn}(t)$.

Since in the decompositions of groups $A$ and $B$ as well, as of
groups $A(t)$ and $B(t)$, into amalgamated free product stated
above the amalgamated subgroups are contained in the centre of
each free factor, by Proposition 1.4 we have

\proclaim{\indent Proposition 2.1} For any element $g$ of group
$A$ (or $B$) not belonging to subgroup $H$ the equality
$g^{-1}Hg\cap H=\langle c\rangle$ (respectively, $g^{-1}Hg\cap
H=\langle d\rangle$) holds. In particular, for any element $g$ of
group $A$ or $B$ and for any element $h$ of subgroup $H$ element
$g^{-1}hg$ belongs to subgroup $H$ if and only if elements $g$ and
$h$ commute.

Similarly, for any element $g$ of group $A(t)$ (or $B(t)$) not
belonging to subgroup $H(t)$ the equality $g^{-1}H(t)g\cap
H(t)=\langle cH^t\rangle$ (respectively, $g^{-1}H(t)g\cap
H(t)=\langle dH^t\rangle$) holds. In particular, for any elements
$g$ of group $A(t)$ or $B(t)$ and $h$ of subgroup $H(t)$, element
$g^{-1}hg$ belongs to subgroup $H(t)$ if and only if elements $g$
and $h$ commute.
\endproclaim

Propositions 1.1, 1.2 and 2.1 lead to the following criterion for
conjugacy of elements of groups $G_{mn}$ and $G_{mn}(t)$:

\proclaim {\indent Proposition 2.2} Let $G$ be any of groups
$G_{mn}$ and $G_{mn}(t)$. Every element of $G$ is conjugate to a
cyclically reduced element. If lengths of two cyclically reduced
elements are nonequal then these elements are not conjugate in
$G$. Let $f$ and $g$ be cyclically reduced elements of $G$ such
that $l(f)=l(g)$. Then \roster \item If $l(f)=l(g)=1$ then $f$ and
$g$ are conjugate in group $G$ if and only if $f$ and $g$ belong
to the same free factor and are conjugate in this factor. \item If
$l(f)=l(g)>1$ then $f$ and $g$ are conjugate in group $G$ if and
only if element $f$ is $H$-conjugate or $H(t)$-conjugate to some
cyclic permutation of $g$.
\endroster
\endproclaim

We now consider some further properties of groups $G_{mn}$ and
$G_{mn}(t)$. The following assertion is easily checked.

\proclaim {\indent Proposition 2.3} For any homomorphism $\varphi$
of group $G_{mn}$ onto a finite group $X$ there exist an integer
$t>1$ and a homomorphism $\psi$ of group $G_{mn}(t)$ onto group
$X$ such that $\varphi=\rho_t\psi$.

Also, for any integers $t>1$ and $s>1$ such that $t$ divides $s$
there exists homomorphism $\varphi:G_{mn}(s) \to G_{mn}(t)$ such
that $\rho_t =\rho_s\varphi$.

The same assertions are valid for groups $A$ and $B$.
\endproclaim

\proclaim {\indent Proposition 2.4} For any element $g$ of group
$A$ or $B$, if $g$ does not belong to subgroup $H$ then there
exists an integer $t_0>1$ such that for any positive integer $t$,
divisible by $t_0$, element $g\rho_t$ does not belong to subgroup
$H\rho_t$.
\endproclaim

\demo{\indent Proof} We shall assume that $g\in A$; the case when
$g\in B$ can be treated similarly.

So, let element $g\in A$ do not belong to subgroup $H$ and let
$g=x_1x_2\cdots x_r$ be a reduced form of $g$ (in the
decomposition of group $A$ into amalgamated free product).

If $r=1$ then element $g$ belongs to subgroup $\langle a\rangle$
i.~e. $g=a^k$ for some integer $k$. Since $g\notin H$, integer $k$
is not divisible by $m$. Then for any integer $t>0$, in group
$A(t)$, element $g\rho_t$ of subgroup $\langle a;\
a^{mt}=1\rangle$ cannot belong to the amalgamated subgroup
(generated by element $a^m$) of the decomposition of group $A(t)$
and consequently cannot belong to the free factor $H(t)=H\rho_t$.

Let $r>1$. Every component $x_i$ of the reduced form of element
$g$ is either of form $a^k$, where integer $k$ is not divisible by
$m$, or of form $c^kd^l$ where $l\neq 0$. If integer $t_0$ is
chosen such that the exponent $l$ of any component $x_i$ of the
second kind is not divisible by $t_0$, then for any $t$ divisible
by $t_0$ the image $x_i\rho_t=x_iH^t$ of such component will not
belong to the amalgamated subgroup of group $A(t)$ (generated,
let's remind, by element $cH^t$). Moreover, as above in case
$r=1$, the images of components of the first kind will not belong
to the amalgamated subgroup. Therefore, the form
$$
g\rho_t = (x_1\rho_t)(x_2\rho_t) \cdots (x_r\rho_t)
$$
of element $g\rho_t$ is reduced in group $A(t)$ and since $r>1$,
$g\rho_t$ does not belong to the free factor $H(t)=H\rho_t$. The
proposition is proven.
\enddemo

Proposition 2.4 obviously implies the

\proclaim {\indent Proposition 2.5} For any element $g$ of group
$G_{mn}$ there exists an integer $t_0>1$ such that for all
positive integer $t$, divisible by $t_0$, the length of element
$g\rho_t$ in group $G_{mn}(t)$ coincides with the length of
element $g$ in group $G_{mn}$.
\endproclaim

\proclaim {\indent Proposition 2.6} For any elements $f$ and $g$
of group $A$ (or $B$) such that element $f$ does not belong to the
double coset $HgH$, then there exists an integer $t_0>1$ such that
for any positive integer $t$, divisible by $t_0$, element
$f\rho_t$ does not belong to the double coset $H(t)(g\rho_t)H(t)$.
\endproclaim

\demo{\indent Proof} We can again consider only the case when
elements $f$ and $g$ belong to subgroup $A$. So, let us suppose
that element $f\in A$ does not belong to the double coset $HgH$.
In view of Proposition 2.4, it is enough to prove that there
exists a homomorphism $\varphi$ of $A$ onto a finite group $X$
such that element $f\varphi$ of $X$ does not belong to the double
coset $(H\varphi)(g\varphi)(H\varphi)$.

To this end let's consider the quotient group $\overline A=A/C$ of
group $A$ by its (central) subgroup $C =\langle c\rangle$. The
image $xC$  of an element $x\in A$ in group $\overline A$ will be
denoted by $\overline x$.

It is obvious that group $\overline A$ is the (ordinary) free
product of cyclic group $X$ of order $m$, generated by element
$\overline a$, and the infinite cycle $Y$, generated by $\overline
d$. The canonical homomorphism of group $A$ onto group $\overline
A$ maps subgroup $H$ onto subgroup $Y$ and consequently, the image
of the double coset $HgH$ is the double coset $Y\overline gY$.
Since $C\leqslant H$, the element $\overline f$ does not belong to
this coset.

We can assume, without loss of generality, that any element
$\overline f$ and $\overline g$, if it is different from identity,
has reduced form the first and the last syllables of which do not
belong to subgroup $Y$. Every $Y$-syllable of these reduced forms
is of form $\overline d^k$ for some integer $k\neq 0$. Since the
set $M$ of all such exponents $k$ is finite, we can choose an
integer $t>0$ such that for any $k\in M$ the inequality $t>2|k|$
holds. Let's denote by $\widetilde A$ the factor group of group
$\overline A$ by the normal closure of element $\overline d^t$.
Group $\widetilde A$ is the free product of groups $X$ and $Y/Y^t$
and since different integers from $M$ are not relatively congruent
and are not congruent to zero modulo $t$, then the reduced forms
of the images $\widetilde f$ and $\widetilde g$ of elements
$\overline f$ and $\overline g$ in group $\widetilde A$ are the
same as in group $\overline A$. In particular, element $\widetilde
f$ does not belong to the double coset $Y/Y^t\,\widetilde
g\,Y/Y^t$. Since this coset consists of a finite number of
elements and group $\widetilde A$ is residually finite, then there
exists a normal subgroup $\widetilde N$ of finite index of group
$\widetilde A$ such, that
$$
\widetilde f\notin (Y/Y^t\widetilde gY/Y^t) \cdot \widetilde N.
$$
If now $\theta$ is the product of the canonical homomorphisms of
group $A$ onto group $\overline A$ and of group $\overline A$ onto
group $\widetilde A$ and $N$ is the full pre-image by $\theta$ of
subgroup $\widetilde N$, then $N$ is a normal subgroup of finite
index of group $A$ and $f\notin (HgH) N$.  Thus, the canonical
homomorphism $\varphi$ of group $A$ onto quotient group $A/N$ has
the required property and the proposition is proven.
\enddemo
\bigskip

\centerline{ \bf 3. Proof of Theorem}
\medskip

We prove first the following proposition.

\proclaim {\indent Proposition 3.1} If elements $f$ and $g$ of
group $G_{mn}$ such that $l(f) =l(g)>1$ are not $H$-conjugate,
then for some integer $t>1$, elements $f\rho_t$ and $g\rho_t$ of
group $G_{mn}(t)$ are not $H(t)$-conjugate.
\endproclaim

\demo{\indent Proof} Let $f=x_1x_2 \cdots x_r$ and $g=y_1y_2\cdots
y_r$ be the reduced forms in group $G_{mn} = (A*B;\ H)$ of
elements $f$ and $g$.

We remind (see Proposition 1.3) that elements $f$ and $g$ are
$H$-conjugate if and only if there exist elements $h_0$, $h_1$,
$h_2$, \dots, $h_r=h_0$ of $H$ such that for any $i=1, 2, \dots,
r$, we have
$$
x_i=h_{i-1}^{-1}y_ih_i.
$$
It then follows, in particular, that for each $i=1, 2, \dots, r$
elements $x_i$ and $y_i$ should lie in the same free factor $A$ or
$B$ and define the same double coset modulo $(H, H)$. So, we
consider separately some cases.

{\it Case 1.} Suppose that for some index $i$ elements $x_i$ and
$y_i$ lie in different free factors $A$ and $B$ of group $G_{mn}$
(and, certainly, are not in subgroup $H$). It follows from
Proposition 2.4 that there exists an integer $t>1$ such that
elements $x_i\rho_t$ and $y_i\rho_t$ do not belong to the same
free factor $A(t)$ or $B(t)$ of group $G_{mn}(t)$ (and, as before,
lie in free factors of this group). Hence, by Proposition 1.3, in
group $G_{mn}(t)$ elements $f\rho_t$ and $g\rho_t$ are not
$H(t)$-conjugate.
\medskip

{\it Case 2.} Let now for any $i=1, 2, \dots, r$ elements $x_i$
and $y_i$ belong to the same subgroup $A$ or $B$  and for some $i$
element $x_i$ does not belong to the double coset $Hy_iH$. By
Proposition 2.6, there exists an integer $t_1>1$ such that for any
positive integer $t$, divisible by $t_1$, element $x_i\rho_t$ is
not in the double coset $H(t)(y_i\rho_t)H(t)$. From Proposition
2.5, there exist integers $t_2>1$ and $t_3>1$ such that for any
positive integer $t$ divisible by $t_2$ the length of element
$f\rho_t$ in group $G_{mn}(t)$ is equal to $r$ and for any
positive integer $t$, divisible by $t_3$, the length of element
$g\tau_t$ in group $G_{mn}(t)$ is equal to $r$. Thus, if
$t=t_1t_2t_3$ then in group $G_{mn}(t)$, elements $f\rho_t$ and
$g\rho_t$ have the reduced forms
$$
(x_1\rho_t)(x_2\rho_t) \cdots (x_r\rho_t) \quad \text{and} \quad
(y_1\rho_t)(y_2\rho_t) \cdots (y_r\rho_t)
$$
respectively and element $x_i\rho_t$ is not in the double coset
$H(t)(y_i\rho_t)H(t)$. Again, by Proposition 1.3 these elements
are not $H(t)$-conjugate.
\medskip

{\it Case 3.} We now consider the case, when for any $i=1, 2,
\dots, r$ elements $x_i$ and $y_i$ lie in the same free factor $A$
or $B$ and also determine the same double coset modulo $(H, H)$.
We prove some lemmas.

\proclaim{\indent Lemma 1} Let elements $x$ and $y$ of one of
groups $A$ or $B$ do not belong to subgroup $H$ and $x\in HyH$,
i.~e.
$$
x=c^{\alpha}d^{\beta}yc^{\gamma}d^{\delta} \tag2
$$
for some integers $\alpha$, $\beta$, $\gamma$ and $\delta$. If
elements $x$ and $y$ belong to group $A$, then integers
$\alpha+\gamma$, $\beta$ and $\delta$ are uniquely determined by
the equality (2). If elements $x$ and $y$ belong to group $B$,
then integers $\beta +\delta$, $\alpha$ and $\gamma$ are uniquely
determined by equality (2).
\endproclaim

\demo{\indent Proof} Let elements $x$ and $y$ belong to subgroup
$A$ and let $ x=c^{\alpha_1}d^{\beta_1}yc^{\gamma_1}d^{\delta_1}$
and $ x=c^{\alpha_2}d^{\beta_2}yc^{\gamma_2}d^{\delta_2}$ for some
integers $\alpha_1$, $\beta_1$, $\gamma_1$, $\delta_1$,
$\alpha_2$, $\beta_2$, $\gamma_2$ and $\delta_2$. Rewriting the
equality $c^{\alpha_1}d^{\beta_1}yc^{\gamma_1}d^{\delta_1} =
 c^{\alpha_2}d^{\beta_2}yc^{\gamma_2}d^{\delta_2}$
as $y^{-1}c^{\alpha_1-\alpha_2}d^{\beta_1-\beta_2}y =
 c^{\gamma_2-\gamma_1}d^{\delta_2-\delta_1}$,
then, by Proposition 2.1, we conclude that
$c^{\alpha_1-\alpha_2}d^{\beta_1-\beta_2} =
 c^{\gamma_2-\gamma_1}d^{\delta_2-\delta_1}
$ and also that this element should belong to subgroup $\langle c
\rangle$, i.~e. $\beta_1-\beta_2 =\delta_2-\delta_1=0$. So, we
have $c^{\alpha_1-\alpha_2} = c^{\gamma_2-\gamma_1}$ and since the
order of element $c$ is infinite, then $\alpha_1-\alpha_2 =
\gamma_2-\gamma_1$.

Thus, $\beta_1 =\beta_2$, $\delta_1 =\delta_2$ and $\alpha_1
+\gamma_1 = \alpha_2 +\gamma_2$ as it was required. The case when
elements $x$ and $y$ belong to group $B$ is esteemed similarly.
\enddemo

\proclaim{\indent Lemma 2} Let elements $x$ and $y$ of one of
groups $A(t)$ or $B(t)$ do not belong to subgroup $H(t)$ and $x\in
H(t)yH(t)$, i.~e.
$$
 x=(cH^t)^{\alpha}(dH^t)^{\beta}y(cH^t)^{\gamma}(dH^t)^{\delta} \tag3
$$
for some integers $\alpha$, $\beta$, $\gamma$ and $\delta$. If
elements $x$ and $y$ belong to group $A(t)$, then integers
$\alpha+\gamma$, $\beta$ and $\delta$ are uniquely determined
modulo $t$ by equality (3). If elements $x$ and $y$ belong to
group $B(t)$, then integers $\beta +\delta$, $\alpha$ and $\gamma$
are uniquely determined modulo $t$ by equality (3).
\endproclaim

The proof of lemma 2 is completely similar to that of lemma 1.

\proclaim{\indent Lemma 3} Let $f=x_1x_2 \cdots x_r$ and
$g=y_1y_2\cdots y_r$ be reduced forms of elements $f$ and $g$ of
group $G_{mn}$, where $r>1$ and let for every $i=1, 2, \dots, r$
the equality $x_i=u_iy_iv_i$ holds, for some elements $u_i$ and
$v_i$ of subgroup $H$. Then there exists at most one sequence
$h_0$, $h_1$, $h_2$, \dots, $h_r$ of elements of subgroup $H$ such
that for any $i=1, 2, \dots, r$
$$
x_i=h_{i-1}^{-1}y_ih_i.\tag4
$$
Moreover, if $u_i=c^{\alpha_i}d^{\beta_i}$ and
$v_i=c^{\gamma_i}d^{\delta_i}$ for some integers $\alpha_i$,
$\beta_i$, $\gamma_i$ and $\delta_i$ ($i=1, 2, \dots, r$) and
$x_1, y_1\in A$, then such sequence exists if, and only if, for
any $i$, $1<i<r$,
$$
\aligned
 \alpha_i+\alpha_{i+1}+\gamma_{i-1}+\gamma_i &= 0, \
\text {if $i$ is odd, and} \\
 \beta_i+\beta_{i+1}+\delta_{i-1}+\delta_i &= 0, \ \text{if $i$ is
even.}
\endaligned \tag5
$$
\endproclaim

\demo{\indent Proof} We suppose first that the sequence $h_0$,
$h_1$, $h_2$, \dots, $h_r$ of elements of subgroup $H$, satisfying
equality (4), exists and let's write its elements as
$h_i=c^{\mu_i}d^{\nu_i}$, for some integers $\mu_i$ and $\nu_i$.

Then, since for any $i=1, 2, \dots, r$ the equality
$h_{i-1}^{-1}y_ih_i=u_iy_iv_i$ holds, we have
$$
 y_i^{-1}(h_{i-1}u_i)y_i=h_iv_i^{-1}.\tag6
$$
As $r>1$, then every element $y_i$, belonging to one of the
subgroups $A$ or $B$, does not lie in subgroup $H$, and
consequently, from Proposition 2.1, for any $i=1, 2, \dots, r$, we
have the equality $h_{i-1}u_i=h_iv_i^{-1}$. Substituting the
expressions of elements $h_i$, $u_i$ and $v_i$, we have for every
$i=1, 2, \dots, r$  $c^{\mu_{i-1}+\alpha_i}d^{\nu_{i-1}+\beta_i}=
c^{\mu_i-\gamma_i} d^{\nu_i-\delta_i}$ and hence we obtain the
system of numeric equations
$$
\mu_{i-1} + \alpha_i =\mu_i-\gamma_i, \quad \nu_{i-1} + \beta_i
=\nu_i-\delta_i \qquad ( i=1, 2, \dots, r). \tag7
$$
Since by (6) for every $i=1, 2, \dots, r$ element $h_iv_i^{-1}$
belongs to the intersection \linebreak $y_i^{-1}Hy_i\cap H$ and,
by supposition, elements $y_1$, $y_3$, \dots belong to subgroup
$A$ and elements $y_2$, $y_4$, \dots belong to subgroup $B$, then
from Proposition 2.1, it follows that for odd $i$, we should have
$h_{i-1} u_i=h_iv_i^{-1} \in \langle c \rangle$, and for even $i$,
we should have $h_{i-1} u_i=h_iv_i^{-1} \in \langle d \rangle$. It
means that for every odd $i$ ($i=1, 2, \dots, r$) we have the
equalities $\nu_{i-1} + \beta_i=0$ and $\nu_i-\delta_i=0$, and for
every even $i$ ($i=1, 2, \dots, r $) we have the equalities
$\mu_{i-1} + \alpha_i=0$ and $\mu_i-\gamma_i=0$.

Hence, the values of the integers $\mu_i$ and $\nu_i$ are
determined uniquely. Namely,
$$
\mu_i =\cases \gamma_i, \ \text{if $i$ is even and
$2\leqslant i \leqslant r$} \\
-\alpha_{i+1}, \ \text{if $i$ is odd and $1\leqslant i \leqslant
r-1 $,}
\endcases \tag8
$$
and
$$
\nu_i =\cases -\beta_{i+1}, \ \text{if $i$ is even and
$0\leqslant i \leqslant r-1$} \\
\delta_i, \ \text{if $i$ is odd and $1\leqslant i \leqslant r $,}
\endcases \tag9
$$
Moreover, from the equalities (7) it follows, that
$$
\aligned
 \mu_0&=-(\alpha_1+\alpha_2+\gamma_1), \ \text{and} \\
 \nu_r&=\beta_r+\delta_{r-1}+\delta_r, \ \text{if $r$ is even,} \\
 \mu_r&=\alpha_r+\gamma_{r-1}+\gamma_r, \ \text{if $r$ is odd.}
\endaligned \tag10
$$

Thus, the statement that there can exist at most one sequence of
elements $h_0$, $h_1$, $h_2$, \dots, $h_r$ of subgroup $H$
satisfying the equalities (4) is demonstrated.

Substituting the value $\mu_{i-1} = \gamma_{i-1}$ and $\mu_i
=-\alpha_{i+1}$ defined in (8) in the equalities $\mu_{i-1} +
\alpha_i =\mu_i-\gamma_i$ of system (7), where $1<i<r$ and $i$ is
odd, we obtain $\alpha_i+\alpha_{i+1}+\gamma_{i-1}+\gamma_i$.
Similarly, substituting the value $\nu_{i-1} = \delta_{i-1}$ and
$\nu_i =-\beta_{i+1}$ defined in (9) in the equalities $\nu_{i-1}
+ \beta_i =\nu_i-\delta_i$ of systems (7), where $1<i<r$ and $i$
is even, we obtain $\beta_i+\beta_{i+1}+\delta_{i-1}+\delta_i =
0$.

Thus, under the existence in subgroup $H$ of sequence $h_0$,
$h_1$, $h_2$, \dots, $h_r$ of elements satisfying (4), conditions
(5) are satisfied.

Conversely, suppose conditions (5) are satisfied. Let
$$
h_0=c^{-(\alpha_1+\alpha_2+\gamma_1)}d^{-\beta_1}
$$
and for all indexes $i$ such that $1\leqslant i < r$, we set
$$
h_i =\cases c^{\gamma_i}d^{-\beta_{i+1}}, &\text{if $i$ is even}, \\
c^{-\alpha_{i+1}}d^{\delta_i}, &\text{if $i$ is odd}.
\endcases
$$
At last, for $i=r$ we set
$$
h_r =\cases c^{\gamma_r}d^{\beta_r+\delta_{r-1}+\delta_r},
 &\text{if $r$ is even}, \\
 c^{\alpha_r+\gamma_{r-1}+\gamma_r}d^{\delta_r}, &\text{if $r$ is odd}.
\endcases
$$

Let's show that, the so defined sequence of elements $h_0$, $h_1$,
$h_2$, \dots, $h_r$ really fits to the equalities (4).

If $i=1$, using the  permutability of elements $c$ and $y_1$ we
have
$$
h_0^{-1}y_1h_1 =
 c^{\alpha_1+\alpha_2+\gamma_1}d^{\beta_1}y_1c^{-\alpha_2}d^{\delta_1}=
 c^{\alpha_1}d^{\beta_1}y_1c^{\gamma_1}d^{\delta_1}=x_1.
$$

If $1<i<r$ and integer $i$ is even, using the equality
$\beta_i+\beta_{i+1}+\delta_{i-1}+\delta_i = 0$ and permutability
of elements $d$ and $y_i$, we have
$$
\multline h_{i-1}^{-1}y_ih_i =
 c^{\alpha_i}d^{-\delta_{i-1}}y_ic^{\gamma_i}d^{-\beta_{i+1}}
 =c^{\alpha_i}y_ic^{\gamma_i}d^{-(\delta_{i-1}+\beta_{i+1})} \\
 =c^{\alpha_i}y_ic^{\gamma_i}d^{\beta_i+\delta_i}
 =c^{\alpha_i}d^{\beta_i}y_ic^{\gamma_i}d^{\delta_i}=x_i.
\endmultline
$$

If $1 < i < r$ and integer $i$ is odd, using the equality
$\alpha_i+\alpha_{i+1}+\gamma_{i-1}+\gamma_i=0$ and permutability
of elements $c$ and $y_i$, we have
$$
\multline h _ {i-1} ^ {-1} y_ih_i =
 c^{-\gamma_{i-1}}d^{\beta_i}y_ic^{-\alpha_{i+1}}d^{\delta_i}
 =c^{-(\gamma_{i-1}+\alpha_{i+1})}d^{\beta_i}y_id^{\delta_i} \\
 =c^{\alpha_i+\gamma_i}d^{\beta_i}y_id^{\delta_i}
 =c^{\alpha_i}d^{\beta_i}y_ic^{\gamma_i}d^{\delta_i}=x_i.
\endmultline
$$

If integer $r$ is even, then
$$
h_{r-1}^{-1}y_rh_r =
c^{\alpha_r}d^{-\delta_{r-1}}y_rc^{\gamma_r}d^{\beta_r+\delta_{r-1}+\delta_r}
 =c^{\alpha_r}d^{\beta_r}y_rc^{\gamma_r}d^{\delta_r}=x_r,
$$
and if $r$ is odd, then
$$
h_{r-1}^{-1}y_rh_r =
c^{-\gamma_{r-1}}d^{\beta_r}y_rc^{\alpha_r+\gamma_{r-1}+\gamma_r}d^{\delta_r}
 =c^{\alpha_r}d^{\beta_r}y_rc^{\gamma_r}d^{\delta_r}=x_r.
$$
Lemma 3 is completely demonstrated.
\enddemo

Similar argument gives the

\proclaim{\indent Lemma 4} Let $f=x_1x_2 \cdots x_r$ and $g=y_1y_2
\cdots y_r$ be reduced forms of elements $f$ and $g$ of group
$G_{mn}(t)$, where $r>1$, and let for every $i=1, 2, \dots, r$ the
equality $x_i=u_iy_iv_i$ takes place, for some elements $u_i$ and
$v_i$ of subgroup $H(t)$. Then there exists at most one sequence
$h_0$, $h_1$, $h_2$, \dots, $h_r$ of elements of subgroup $H(t)$
such that for any $i=1, 2, \dots, r$, we have
$$
x_i=h_{i-1}^{-1}y_ih_i.\tag11
$$
Moreover, if $u_i = (cH^t)^{\alpha_i}(dH^t)^{\beta_i}$ and $v_i=
(cH^t)^{\gamma_i}(dH^t)^{\delta_i}$ for some integers $\alpha_i$,
$\beta_i$, $\gamma_i$ and $\delta_i$ ($i=1, 2, \dots, r$) and
$x_1, y_1\in A(t)$, then such sequence exist if and only if for
any integer $i$, $1<i<r$, we have
$$
\aligned
 \alpha_i+\alpha_{i+1}+\gamma_{i-1}+\gamma_i &\equiv 0 \pmod t, \
\text{if $i$ is odd, and} \\
 \beta_i+\beta_{i+1}+\delta_{i-1}+\delta_i &\equiv 0 \pmod t,
\ \text{if $i$ is even.}
\endaligned \tag12
$$
\endproclaim

We can now end the consideration of case 3 and thus complete the
proof of Proposition 3.1.

By Proposition 1.3, elements $f$ and $g$ are $H$-conjugate if and
only if their cyclic permutation $x_2 \cdots x_rx_1$ and $y_2
\cdots y_ry_1$ are $H$-conjugate. Therefore, we can suppose,
without loss of generality, that elements $x_1$ and $y_1$ belong
to subgroup $A$.

By supposition, for every $i=1, 2, \dots, r$ there exist integers
$\alpha_i$, $\beta_i$, $\gamma_i$ and $\delta_i$ such that
$$
 x_i=c^{\alpha_i}d^{\beta_i}y_ic^{\gamma_i}d^{\delta_i}.
$$

If  integer $i$, $1<i<r$, is even, then $y_i\in B$, $y_{i-1},
y_{i+1} \in A$ and consequently, from Lemma 1, integers $\beta_i
+\delta_i$, $\delta_{i-1}$ and $\beta_{i+1}$ do not depend from
the particular expressions of the elements $x_i$ (in the form
$x_i=uy_iv$, where $u, v\in H$), and these integers are uniquely
determined by the sequences $x_1$, $x_2$, \dots, $x_r$ and $y_1$,
$y_2$, \dots, $y_r$. Similarly, for any odd $i$, $1<i<r$, integers
$\alpha_i +\gamma_i$, $\gamma_{i-1}$ and $\alpha_{i+1}$ are
uniquely determined by these sequences. It means, in turn, that
the satisfiability of conditions (5) of lemma 3 depends only on
the sequences $x_1$, $x_2$, \dots, $x_r$ and $y_1$, $y_2$, \dots,
$y_r$.

Let's now conditions (5) of Lemma 3 are not satisfied, i.~e.
either for some even $i$, $1<i<r$, the sum
$\beta_i+\beta_{i+1}+\delta_{i-1}+\delta_i$ is different from
zero, or for some odd $i$, $1<i<r$, the sum
$\alpha_i+\alpha_{i+1}+\gamma_{i-1}+\gamma_i$ is different from
zero. Then it is possible to find an integer $t_1>1$, not dividing
the respective sum. Let's also choose, according to Proposition
2.5, the integers $t_2 > 1$ and $t_3 > 1$ such that for all
positive integer $t$, divisible by $t_2$, the length of element
$f\rho_t$ in group $G_{mn}(t)$ is equal to $r$ and for all
positive integer $t$, divisible by $t_3$, the length of element
$g\rho_t$ in group $G_{mn}(t)$ is equal to $r$. Then if
$t=t_1t_2t_3$, in group $G_{mn}(t)$, elements $f\rho_t$ and
$g\rho_t$ have reduced forms
$$
(x_1\rho_t)(x_2\rho_t) \cdots (x_r\rho_t) \quad \text{and} \quad
(y_1\rho_t)(y_2\rho_t) \cdots (y_r\rho_t)
$$
respectively. Further, for every $i=1, 2, \dots, r$,
$$
 x_i\rho_t=(cH^t)^{\alpha_i}(dH^t)^{\beta_i}(y_i\rho_t)
(cH^t)^{\gamma_i}(dH^t)^{\delta_i}.
$$
From Lemma 2 and the selection of integer $t_1$ if follows that
for elements $f\rho_t$ and $g\rho_t$, conditions (12) of Lemma 4
are not satisfied. Therefore from this lemma and Proposition 1.3,
elements $f\rho_t$ and $g\rho_t$ are not $H(t)$-conjugate in group
$G_{mn}(t)$.

Let now the reduced forms of elements $f$ and $g$ satisfy
conditions (5) of Lemma 3. Then according to this lemma, in
subgroup $H$ there exists the only sequence of elements $h_0$,
$h_1$, $h_2$, \dots, $h_r$ such that for any $i=1, 2, \dots, r$,
we have $x_i=h_{i-1}^{-1}y_ih_i$. Since elements $f$ and $g$ are
not $H$-conjugate, then elements $h_0$ and $h_r$ should be
different. Since group $G_{mn}$ is residually finite, by
Proposition 2.3 there exists an integer $t_1>1$ such that
$h_0\rho_{t_1} \neq h_r\rho_{t_1}$. Let's choose one more integer
$t_2$ such that in group $G_{mn}(t_2)$ elements $f\rho_{t_2}$ and
$g\rho_{t_2}$ have length $r$. Then if $t=t_1t_2$, in group
$G_{mn}(t)$, $h_0\rho_{t} \neq h_r\rho_{t}$, elements $f\rho_{t}$
and $g\rho_{t}$ have reduced forms
$$
(x_1\rho_t)(x_2\rho_t) \cdots (x_r\rho_t) \quad \text{and} \quad
(y_1\rho_t)(y_2\rho_t) \cdots (y_r\rho_t)
$$
respectively and for every $i=1, 2, \dots, r$
$$
 x_i\rho_t=(cH^t)^{\alpha_i}(dH^t)^{\beta_i}(y_i\rho_t)
(cH^t)^{\gamma_i}(dH^t)^{\delta_i}.
$$
Moreover, in group $G_{mn}(t)$, for any $i=1, 2, \dots, r$, the
equality
$$
 x_i\rho_t=(h_{i-1}\rho_t)^{-1}(y_i\rho_t)(h_i\rho_t)
$$
holds.

Since the sequence $h_0\rho_t$, $h_1\rho_t$, $h_2\rho_t$, \dots,
$h_r\rho_t$ is, by Lemma 4, the unique sequence of elements of
subgroup $H(t)$ satisfying these equalities, then from Proposition
1.3, elements $f\tau_t$ and $g\tau_t$ are not $H(t)$-conjugate in
group $G_{mn}(t)$.

Hence Proposition 3.1 is completely demonstrated.
\enddemo

We now proceed directly to the proof of Theorem.

As remarked above, for any $t>1$, group $G_{mn}(t)$ is conjugacy
separable. Therefore, for the proof of the Theorem it is enough to
show that for any two non-conjugate in group $G_{mn}$ elements $f$
and $g$ of group $G_{mn}$, there exists an integer $t>1$ such that
elements $f\rho_t$ and $g\rho_t$ are not conjugate in group
$G_{mn}(t)$. To this end we make use of criterion of conjugacy of
elements of groups $G_{mn}$ and $G_{mn}(t)$ given in Proposition
2.2.

So, let $f=x_1x_2 \cdots x_r$ and $g=y_1y_2 \cdots y_r$ be the
reduced forms in group $G_{mn} = (A*B;\ H)$ of two non-conjugate
in group $G_{mn}$ elements $f$ and $g$. By Proposition 2.2 we can
assume that elements $f$ and $g$ are cyclically reduced. In
accordance with this proposition we must consider separately some
cases.
\smallskip

1. $l(f)\neq l(g)$.

From Proposition 2.5, there exists an integer $t_1>1$ such that
for any positive integer $t$, divisible by $t_1$, the length of
element $f\rho_t$ in group $G_{mn}(t)$ coincides with the length
of element $f$ in group $G_{mn}$. Similarly, there exists integer
$t_2>1$ such that for any positive integer $t$, divisible by
$t_2$, the length of element $g\rho_t$ in group $G_{mn}(t)$
coincides with the length of element $g$ in group $G_{mn}$. Thus,
if $t=t_1t_2$ elements $f\rho_t$ and $g\rho_t$ of group
$G_{mn}(t)$ are, as it is easy to see, cyclically reduced and have
different length. Hence, by Proposition 2.2, these elements are
not conjugate in this group. Thus integer $t$ is the required.
\smallskip

2. $l(f) = l(g) = 1$.

In this case each of these elements belongs to one of the
subgroups $A$ or $B$. Suppose first that both elements lie in the
same of these subgroups; let it be, for instance, subgroup $A$.
Since elements $f$ and $g$ are not conjugate in group $A$ and
group $A$ is conjugacy separable [4], then by Proposition 2.3
there exists an integer $t>1$ such that in group $A(t)$ elements
$f\rho_t$ and $g\rho_t$ are not conjugate. Proposition 2.2 now
implies that elements $f\rho_t$ and $g\rho_t$ are not conjugate in
group $G_{mn}(t)$.

Let now element $f$ belongs to subgroup $A$, element $g$ belongs
to subgroup $B$ and these elements don't belong to subgroup $H$.
By Proposition 2.4, there exist integers $t_1>1$ and $t_2>1$ such
that for any positive integer $t$, divisible by $t_1$, element
$f\rho_t$ does not belong to subgroup $H\rho_t$, and for any
positive integer $t$, divisible by $t_2$, element $g\rho_t$ does
not belong to subgroup $H\rho_t$. Then if $t=t_1t_2$, Proposition
2.2 implies that elements $f\rho_t$ and $g\rho_t$ are not
conjugate in group $G_{mn}(t)$.
\smallskip

3. $l(f) =l(g)>1$.

Let $g_i=y_iy_{i+1} \cdots y_ry_1 \cdots y_{i-1}$ ($i=1, 2, \dots,
r$) be all the cyclic permutations of element $g$. Since elements
$g$ and $g_i$ are conjugate and elements $f$ and $g$ are not
conjugate, element $f$ is not $H$-conjugate to any of elements
$g_1$, $g_2$, \dots, $g_r$. It follows from Proposition 3.1 that
for every $i=1, 2, \dots, r$, there exists an integer $t_i>1$ such
that elements $f\rho_{t_i}$ and $g_i\rho_{t_i}$ are not
$H(t_i)$-conjugate in group $G_{mn}(t_i)$. Let also integer $t_0
> 1$ be chosen such that for all positive integer $t$, divisible by
$t_0$, elements $f\rho_{t}$ and $g\rho_{t}$ have length $r$ in
group $G_{mn}(t)$. Then if $t=t_0t_1 \cdots t_r$ in group
$G_{mn}(t)$, elements $f\rho_{t}$ and $g_i\rho_{t}$ have reduced
forms
$$
(x_1\rho_t)(x_2\rho_t) \cdots (x_r\rho_t)
$$
and
$$
(y_i\rho_t)(y_{i+1}\rho_t) \cdots (y_r\rho_t)(y_1\rho_t) \cdots
(y_{i-1}\rho_t)
$$
respectively. Furthermore, elements $f\rho_{t}$ and $g_i\rho_{t}$
are not $H(t)$-conjugate in group $G_{mn}(t)$. Since an arbitrary
cyclic permutation of element $g\rho_{t}$ coincides with some
element $g_i\rho_{t}$, then from Proposition 2.2 it follows that
elements $f\rho_t$ and $g\rho_t$ are not conjugate in group
$G_{mn}(t)$.

The proof of Theorem is now complete.
\end